\documentclass{elsarticle-ppn}

\usepackage{lineno,hyperref}
\modulolinenumbers[5]

\usepackage{textcomp}
\usepackage{amsmath,amsfonts,amsthm,amssymb,amsxtra,bbm}
\usepackage{MnSymbol}
\usepackage{color}

\newtheorem{theorem}{Theorem}
\newtheorem{proposition}[theorem]{Proposition}
\newtheorem{lemma}[theorem]{Lemma}
\newtheorem{corollary}[theorem]{Corollary}

\newcommand{\R}{{\mathbb R}}
\newcommand{\N}{{\mathbb N}}

\newcommand{\be}[1]{\begin{equation}\label{#1}}
\newcommand{\ee}{\end{equation}}

\newcommand{\kin}{{\mathrm{kin}}}

\begin{document}
\begin{frontmatter}

\title{Time averages for kinetic Fokker-Planck equations}

\author{Giovanni M. Brigati\fnref{myfootnote}$^{a,b,}$}
\fntext[myfootnote]{\noindent Email: brigati@ceremade.dauphine.fr \\ Web: \url{https://www.ceremade.dauphine.fr/~brigati/}}



\address[mymainaddress]{CEREMADE, CNRS, UMR 7534, Universit\'e Paris-Dauphine, PSL University, Place du Marechal de Lattre de Tassigny, 75016 Paris, France}
\address[mysecondaryaddress]{Dipartimento di Matematica ``F. Casorati'', Università degli Studi di Pavia, Via Ferrata 5, 27100 Pavia, Italia}

\begin{abstract}
We consider kinetic Fokker-Planck (or Vlasov-Fokker-Planck) equations on the torus with Maxwellian or fat tail local equilibria. Results based on weak norms have recently been achieved by S.~Armstrong and J.-C.~Mourrat in case of Maxwellian local equilibria. Using adapted Poincar\'e and Lions-type inequalities, we develop an explicit and constructive method for estimating the decay rate of time averages of norms of the solutions, which covers various regimes corresponding to subexponential, exponential and superexponential (including Maxwellian) local equilibria. As a consequence, we also derive hypocoercivity estimates, which are compared to similar results obtained by other techniques.
\end{abstract}


\end{frontmatter}

\section{Introduction}\label{Sec1}
Let us consider the \emph{kinetic Fokker-Planck equation}
\begin{equation}
 \label{eq1}
 \partial_t f + v \cdot \nabla_x f = \nabla_v \cdot \left(\nabla_v f + \alpha\,\langle v \rangle^{\alpha-2}\,v\,f \right), \quad \quad f(0,\cdot,\cdot) =f_0. 
 \end{equation}
where $f$ is a function of time $t \geq 0,$ position $x$, velocity $v,$ and $\alpha$ is a positive parameter. Here we use the notation 
$$\langle v \rangle = \sqrt{1 + |v|^2}, \quad \forall\,v\in \mathbb R^d.$$
We consider the spatial domain $\mathbb T:= (0,L)^d \ni x,$ with \emph{periodic boundary conditions}, and define $\Omega_t:=(t,t+\tau) \times \mathbb T,$ for some $\tau >0$, $t\ge0$ and $\Omega=\Omega_0.$ The normalized \emph{local equilibrium}, that is, the equilibrium of the spatially homogeneous case, is
$$\gamma_\alpha(v) = \frac{1}{Z_\alpha}\,\mathrm{e}^{-\langle v \rangle^\alpha}, \quad \forall\,v\in \mathbb R^d,$$
where $Z_\alpha$ is a non-negative normalization factor, so that $d\gamma_\alpha:=\gamma_\alpha(v)\,dv $ is a probability measure. We shall distinguish a sublinear regime if $\alpha \in (0,1),$ a linear regime if $\alpha=1$ and a superlinear regime if $\alpha \geq 1.$ The superlinear regime covers the Maxwellian case $\alpha =2.$ The threshold case $\alpha = 1$ corresponds to a linear growth of $\langle v \rangle^{\alpha}$ as $|v| \to+\infty.$ The estimates in the linear case are similar to the ones of the superlinear regime. In the literature, $\gamma_\alpha$ is said to be \emph{subexponential}, \emph{exponential} or \emph{superexponential} depending whether the regime is sublinear, linear or superlinear.

The mass
$$M:=\iint_{\mathbb T \times \mathbb R^d} f(\cdot,x,v)\,dx\,dv$$
is conserved under the evolution according to the kinetic Fokker-Planck equation~\eqref{eq1}. We are interested in the convergence of the solution to the stationary solution $M\,L^{-d}\,\gamma_\alpha.$ By linearity, we can assume from now on that $M=0$ with no loss of generality. The function 
$$h = \frac{f}{\gamma_\alpha}$$ 
solves the \emph{kinetic-Ornstein-Uhlenbeck} equation
\begin{equation}
 \label{eq2}
 \partial_t h + v \cdot \nabla_x h = \Delta_\alpha h, \quad \quad h(0,\cdot,\cdot) = h_0,
\end{equation}
with $$\Delta_\alpha h := \Delta_v h - \alpha\,v\,\langle v \rangle^{\alpha-2} \cdot \nabla_v h,$$
and \emph{zero-average initial datum} in the sense that
$$\iint_{\mathbb T \times \mathbb R^d}h_0(x,v)\,dx\,d\gamma_\alpha=0.$$
By mass conservation, solutions to ~\eqref{eq2} are zero-average for any $t>0.$

Therefore, we consider the \emph{time average} defined as
$$\strokedint_t^{t+\tau} g(s)\,ds := \frac{1}{\tau} \int_t^{t+\tau} g(s)\,ds$$
without specifying the $\tau$ dependence when not necessary. Our first result is devoted to the decay rate of $h(t,\cdot,\cdot) \to 0$ as $t \to \infty$ using time averages.
\begin{theorem}\label{thm1}
Let $\alpha \geq 1.$ Then, for all $L>0$ and $\tau>0$, there exists a constant $\lambda>0$ such that, for all $h_0 \in \mathrm{L}^2(dx\,d\gamma_\alpha)$ with zero-average, the solution to~\eqref{eq2} satisfies
\begin{equation}
 \label{eq5}
 \strokedint_t^{t+\tau} \|h(s,\cdot,\cdot) \|^2_{\mathrm{L}^2(dx\,d\gamma_\alpha)}\,ds \leq\,\|h_0\|^2_{\mathrm{L}^2(dx\,d\gamma_\alpha)}\,e^{-\lambda\,t}, \quad \forall\,t \geq 0.
\end{equation}
\end{theorem}
The expression of $\lambda$ as a function of $\tau$ and $L$ is given in Section~\ref{Sec4}. To deal with large-time asymptotics in kinetic equations, it is by now standard to use hypocoercivity methods. Although not being exactly a hypocoercive method in the usual sense, Theorem~\ref{thm1} provides us with a hypocoercivity estimate.
\begin{corollary}\label{cor1}
Under the assumptions of Theorem~\ref{thm1}, there exists an explicit constant $C>1$ such that all solutions $h$ to~\eqref{eq2} fulfill 
\begin{equation}\label{hypo}
 \|h(t,\cdot,\cdot)\|^2_{\mathrm{L^2}(dx\,d\gamma_\alpha)} \leq C\,\|h_0\|^2_{\mathrm{L^2}(dx\,d\gamma_\alpha)}\,\mathrm{e}^{-\lambda\,t}, \quad \forall\,t \geq 0.
\end{equation}
\end{corollary}
A typical feature of hypocoercive estimates is the factor $C>1$ in~\eqref{hypo}. 
The prefactor $C>1$ cannot be avoided. Otherwise, inequality \eqref{hypo} would be equivalent to a Poincar\'e inequality where the $\mathrm{L}^2-$norm of a function is controlled with the velocity gradient only. 
We can see explicitly that for $\alpha=2$ the Green function of \eqref{eq1}, computed in \cite{kolmogoroff1934zufallige}, has a built-in delay.
In particular, 
$$\|h_0\|^2_{\mathrm{L}^2(dx\,d\gamma_2)} - \|h(t,\cdot,\cdot)\|^2_{\mathrm{L}^2(dx\,d\gamma_2)} = O(t^3),$$
as $t \to 0^+$.
Note that there is no such a constant in~\eqref{eq5}.

Now, let us turn our attention to the subexponential case $0<\alpha<1$. 
\begin{theorem}\label{thm2}
Let $\alpha \in (0,1)$ Then, for all $L>0$ and $\tau>0$, for all $\sigma >0$, there is a constant $K>0$ such that all solutions to~\eqref{eq2} decay according to
\begin{equation}\label{decay}
\strokedint_t^{t+\tau} \|h(s,\cdot,\cdot)\|^2_{\mathrm{L}^2(dx\,d\gamma_\alpha)}\,ds \leq\,K\,(1+t)^{-\frac\sigma{2\,(1-\alpha)}}\,\iint_{\mathbb T \times \mathbb{R}^d} \langle v\rangle ^\sigma h^2_0\,dx\,d\gamma_\alpha, \quad \forall\,t\geq 0.
\end{equation}
\end{theorem}
Further details will be given in Section~\ref{Sec5}. The constants $\lambda$, $C$ in Corollary~\ref{cor1} and $K$ in Theorem~\ref{thm2} depend on $L>0$ and $\tau>0$ and their values are discussed~later. The rate of Theorem~\ref{thm2} is the same as in the spatially-homogeneous case of~\cite[Proposition~11]{bouin2017hypocoercivity}. In the spatially-inhomogenous case, rates are known, see \cite{MR4265692,bouin2019hypocoercivity}. Finally, see Section~\ref{Sec5} for a discussion of the limit $\alpha\to1^-.$

\medskip Equation~\eqref{eq2} is used in physics to describe the distribution function of a system of particles interacting randomly with some background, see for instance~\cite{balian2007microphysics}. The \emph{kinetic Fokker-Planck equation} is the Kolmogorov forward equation of Langevin dynamics 
\begin{equation*} 
\begin{cases}
dx_t = v_t\,dt, \\ 
dv_t = -v_t + \sqrt{2}\,dW_t,
\end{cases}
\end{equation*}
where $W_t$ is a standard Brownian motion. See~\cite[Introduction]{bernard2020hypocoercivity} for further details on connections with probability theory. The kinetic Fokker-Planck equation~\eqref{eq1} is a simple  kinetic equation which has a long history in mathematics that we will not retrace in details here. Mathematical results go back at least to~\cite{kolmogoroff1934zufallige} and are at the basis of the theory of L.~H\"ormander (see, \emph{e.g.},~\cite{hormander1967hypoelliptic}), at least in the case $\alpha=2.$ For the derivation of the kinetic-Fokker-Planck equation from underlying stochastic ODEs, particularly in the context of astrophysics, we can refer to \cite[eq.~(328)]{RevModPhys.15.1}. Modern hypoellipticity theory emerged from~\cite{MR2034753,MR1787105} and was built up in a fully developed theory in~\cite{villani2009hypocoercivity} with important contributions in~\cite{MR2215889,Mouhot_2006}. Existence theory for solutions to the Vlasov-Fokker-Planck equation was discussed also in \cite[Appendix A]{degond1986global}.

 The word \emph{hypocoercivity} was coined by T.~Gallay, in analogy with the already quoted \emph{hypoelliptic} theory of H\"ormander in~\cite{hormander1967hypoelliptic}. In~\cite{villani2009hypocoercivity}, C.~Villani distinguishes the regularity point of view for elliptic and parabolic problems driven by degenerate elliptic operators from the issue of the long-time behaviour of solutions, which is nowadays attached to the word hypocoercivity. The underlying idea is to twist the reference norm, in order to carry properties (as the coercivity of the operator driving \eqref{eq2}) from velocity direction to space directions, thanks to commutators. Twisting the $\mathrm{H}^1$-norm creates equivalent norms, which are exponentially decaying along the evolution. So works the $\mathrm{H}^1$ framework, see~\cite{villani2009hypocoercivity,villani2006hypocoercive,dolbeault2018varphi}. The $\mathrm{H}^1$-framework has been connected to the \emph{carr\'e du champ} method of D.~Bakry and M.~Emery in~\cite{bakry1985diffusions} by F.~Baudoin, who proved decay also w.r.t.~the Wasserstein distance, as shown in~\cite{baudoin2016wasserstein,baudoin2013bakryp,baudoin2019gamma}. We report also the works \cite{eberle2019couplings,dietert2015convergence}, where accurate convergence rates in the Wasserstein distance for \eqref{eq2} are computed trough a coupling argument.

The $\mathrm{H}^1$ hypocoercivity implies a decay rate for the $\mathrm{L}^2$ norm \cite{Mouhot_2006}, but the corresponding estimates turn out to be sub-optimal. Moreover, kinetic equations driven by non-regularising operators are not well suited for the $\mathrm{H}^1-$framework. This motivates the development of direct $\mathrm{L}^2$ techniques based on a perturbation of the $\mathrm{L}^2$ norm. Such an approach can be found in~\cite{dolbeault2015hypocoercivity} and~\cite{bouin2017hypocoercivity}, which is consistent with diffusion limits. In~\cite{bouin2019hypocoercivity}, the authors extend the technique to the subexponential case. Another possibility is to perform rotations in the phase space and use a Lyapunov inequality for matrices as in~\cite{arnold2014sharp}. This approach gives optimal rates, but it is less general as it requires further algebraic properties for the diffusion operator and a detailed knowledge of its spectrum. The core of~\cite{arnold2014sharp} is a spectral decomposition, that was originally understood via a toy model exposed in~\cite{dolbeault2015hypocoercivity}. In a domain with periodic boundary conditions and no confining potential, the problem is reduced to an infinite set of ODEs corresponding to spatial modes. See~\cite{achleitner2016linear,achleitner2017multi,arnold2021} for details and extensions. 
Other techniques related to hypocoercivity -- involving time-integrated functionals and the application of the so called \textit{kinetic-fluid decomposition}, appear in \cite{strain2004stability,guo2002landau} and subsequent papers. 

A new hypocoercivity theory, involving Sobolev norms with negative exponents of the transport operator, was recently proposed by S.~Armstrong and J.-C.~Mourrat in~\cite{albritton2019variational}. Using space-time adapted Poincar\'e inequalities they derive qualitative hypocoercive estimates in the case $\alpha =2$ on bounded spatial domains.
The constants appearing there are not quantified. One of the difficulties lies in controlling the constant in Lions' Lemma, which is done in our Section~\ref{Sec2}. An extension to the whole space in presence of a confining potential can be found in~\cite{cao2019explicit}. 
Note that the strategy of using time-integrated functionals of the solutions to kinetic equations is present also in \cite{tran2013convergence,duan2011hypocoercivity}.

\medskip Adopting the strategy of \cite{albritton2019variational}, in this paper we study the convergence to equilibrium of solutions to \eqref{eq1} and \eqref{eq2}, as it is a simple benchmark in kinetic theory, \cite{dolbeault2015hypocoercivity, villani2006hypocoercive, albritton2019variational}, 
and a simplified model of the Boltzmann equation when collisions become grazing, see \cite{degond1992fokker}.

Our original contribution lies in making the strategy of \cite{albritton2019variational} effective, and to generalise it to kinetic Fokker-Planck equations where local equilibria are not necessarily Maxwellians. First, we are able to track the Lions' constant in terms of the parameters (see Lemma~\ref{bestlions}). 
Moreover, we achieve a fully constructive proof of the averaging Lemma \ref{lemmaavg}. This allows both for an explicit estimate of the constant and for an adaptation to more general models. One important point is the control in terms of the offset of the solution from the velocity average, without explicitly using gradients, see Proposition \ref{Pro:pgPI}. 
So, we compute explicit and accurate decay rates of time averages of solutions to \eqref{eq2}. Hypocoercivity estimates are obtained as a consequence of these decay rates, see Corollary \ref{cor1}. We perform an analysis for all positive values of~$\alpha,$ which is consistent in the threshold case $\alpha = 1.$ 
Since the estimates are explicit, we are able to compare the strategy of \cite{albritton2019variational} to other $\mathrm{L}^2-$hypocoercivity methods. 

This document is organized as follows. In Section~\ref{Sec2} we collect some preliminary results: Poincar\'e and weighted Poincar\'e inequalities (Propositions~\ref{prop22} and~\ref{weightp}), adapted Lions' inequality (Lemmas~\ref{lemma2} and~\ref{bestlions}). In Section~\ref{Sec3} we introduce an averaging lemma (Lemma~\ref{lemmaavg}), which is then used to prove the generalized Poincar\'e inequality of Proposition~\ref{Pro:pgPI}, at the core of the method. In Section~\ref{Sec4} we use Proposition~\ref{Pro:pgPI} and a Gr\"onwall estimate to prove Theorem~\ref{thm1} and compute an explicit formula for $\lambda$ (Proposition
~\ref{thm1vero}). Section~\ref{Sec5} is devoted to the proof of Theorem~\ref{thm2}, with additional details, and to the limit $\alpha \to 1^-.$ Finally, in Section~\ref{Sec7}, we derive the hypocoercive estimates of Corollary~\ref{cor1}. On the benchmark case $\alpha=2$ in one spatial dimension, we also compare our results with those obtained by more standard methods.

\section{Preliminaries}\label{Sec2}
Let us start with some preliminary results.

\subsection{Weighted spaces}
For functions $g$ of the variable $v$ only, that is, of the so-called homogeneous case, we define the \emph{weighted Lebesgue} and \emph{Sobolev spaces}
$$\mathrm{L}^2_\alpha : = \mathrm{L^2}(\mathbb{R}^d,d\gamma_\alpha)\quad\mbox{and}\quad\mathrm{H}^1_\alpha := \left\{g \in \mathrm{L}^2_\alpha\,:\,\nabla_v g \in \left(\mathrm{L}^2_\alpha\right)^d \right\}.$$
We equip $\mathrm{L}^2_\alpha$ with the scalar product
\be{scal}(g_1,g_2) = \int_{\mathbb{R}^d} g_1(v)\,g_2(v)\,d\gamma_\alpha\ee
and consider on $\mathrm{H}^1_\alpha$ the norm defined by
$$\|h\|^2_{\mathrm{H}^1_\alpha} := \left(\int_{\mathbb R^d} h\,d\gamma_\alpha\right)^2 + \int_{\mathbb R^d} |\nabla_v h|^2\,d\gamma_\alpha$$
as in~\cite{albritton2019variational}. The duality product between $\mathrm{H}^{-1}_\alpha \ni z$ and $\mathrm{H}^1_\alpha \ni g$ is given by
$$\langle z,g \rangle:= \int_{\mathbb{R}^d} \nabla_v w_z \cdot \nabla_v g\,d\gamma_\alpha,$$
where $w_z$ is the weak solution in $\mathrm H^1_\alpha$ to
\begin{equation*}
 -\Delta_\alpha w_z = z - \int_{\mathbb{R}^d} z\,d\gamma_\alpha, \quad \int_{\mathbb{R}^d} w\,d\gamma_\alpha = 0. 
\end{equation*}
Here we write $\int_{\mathbb{R}^d} z\,d\gamma_\alpha$ for functions which are integrable w.r.t.~$d\gamma_\alpha$ and, up to a little abuse of notations, this quantity has to be understood in the distribution sense for more general measures. As a consequence and with the above notations, we define 
$$ \|z\|^2_{\mathrm{H}^{-1}_\alpha} := \left(\int_{\mathbb{R}^d} z\,d\gamma_\alpha\right)^2 + \|w_z\|^2_{\mathrm{H^1_\alpha}}.$$
With these notation, the key property of the operator $\Delta_\alpha,$ is
\begin{equation*}
 \langle g_1 ,\Delta_\alpha g_2 \rangle = -\int \nabla_v g_1 \cdot \nabla_v g_2\,d\gamma_\alpha
\end{equation*}
for any functions $g_1$, $g_2\in\mathrm{H}^1_\alpha.$

\medskip We recall that $\Omega_t=(t,t+\tau) \times \mathbb T\subset\R^+_t\times\R^d_x$ and that $x$-periodic boundary conditions are assumed.  Consider next functions $h$ of $(t,x,v)\in\R^+\times \mathbb T\times\R^d$ and define the space
$$\mathrm{H}_\kin := \left\{ h \in \mathrm{L}^2\left((t,t+\tau)\times \mathbb T; \mathrm{H}^1_\alpha\right)\,:\,\partial_t h + v \cdot \nabla_x h \in \mathrm{L}^2\left(\Omega_t;\mathrm{H}^{-1}_\alpha\right)\,\forall\,t\ge0 \right\}.$$
The dependence of the space on $t,\tau$ is implicit for readability purposes. 
 We can equip $\mathrm{H}_\kin$ with the norm
$$\|h\|^2_\kin := \|h\|^2_{\mathrm{L}^2(\Omega_t;\mathrm{L^2_\alpha}) } + |h|_\kin^2$$
where the \emph{kinetic semi-norm} is given by
$$|h|^2_\kin := \|\nabla_v h\|_{\mathrm{L}^2(\Omega_t;\mathrm{L^2_\alpha}) }^2 + \|\partial_t h + v \cdot \nabla_x h\|_{\mathrm{L}^2(\Omega_t;\mathrm{H^{-1}_\alpha}) }^2.$$ We refer to~\cite[Section 6]{albritton2019variational} for the proof of following result.
\begin{proposition}
The embedding $\mathrm{H}^1_\kin \hookrightarrow \mathrm{L}^2\left(\Omega_t\times \mathbb T; \mathrm L^2_\alpha\right)$ is continuous and compact for any $t\ge0.$
\end{proposition}

\subsection{Poincar\'e inequalities}
In this subsection, we consider functions $g$ depending only on the variable~$v.$
Let $\alpha \geq 1.$ We can state some \emph{Poincar\'e inequalities.}
\begin{proposition}\label{prop22}
If $\alpha\ge1$, there exists a constant $P_\alpha >0$ such that, for all functions $g \in \mathrm{H}^1_\alpha$, we have 
\begin{equation}
 \label{poi}
 \int_{\mathbb R^d}|g-\rho_g|^2\,d\gamma_\alpha \leq P_\alpha \int_{\mathbb R^d} |\nabla_v g|^2\,d\gamma_\alpha\quad\mbox{with}\quad\rho_g:=\int_{\mathbb R^d} g\,d\gamma_\alpha.
\end{equation}
\end{proposition}
With $\alpha\ge1$, the operator $\Delta_\alpha$ admits a compact resolvent on $\mathrm{L}^2(d\gamma_\alpha).$ Then,~\eqref{poi} holds by the standard results of~\cite[Chapter~6]{brezis1973operateurs}. The best constant $P_\alpha$ is such that $P_\alpha^{-1}$ is the minimal positive eigenvalue of~$-\Delta_\alpha.$  See~\cite{cattiaux2019entropic} and the references quoted therein for estimates on $P_\alpha.$ In the case of the Gaussian Poincar\'e inequality, it is shown in~\cite{nash1958continuity} that $P_2 = 1,$ although the result was probably known before.

\subsection{Weighted Poincar\'e inequalities}
Here we consider again functions depending only on $v.$
For $\alpha \in (0,1)$, inequality~\eqref{poi} has to be replaced by the following \emph{weighted Poincar\'e inequality.}
\begin{proposition}\label{weightp}
If $\alpha \in (0,1)$, there exists a constant $P_\alpha>0$ such that, for all functions $g \in \mathrm{H}^1_\alpha$, we have 
\begin{equation}\label{weightpeq}
\int_{\mathbb R^d} \langle v \rangle^{2\,(\alpha-1)}\,|g-\rho_g|^2\,d\gamma_\alpha \leq P_\alpha \int_{\mathbb R^d} |\nabla_v g|^2\,d\gamma_\alpha\quad\mbox{with}\quad\rho_g:=\int_{\mathbb R^d} g\,d\gamma_\alpha. 
\end{equation}
\end{proposition}
For more details, we refer for instance to~\cite[Appendix A]{bouin2019hypocoercivity}. Notice that the average in the l.h.s.~is taken w.r.t.~$d\gamma_\alpha$, not w.r.t.~$\langle v \rangle^{2\,(\alpha-1)}\,d\gamma_\alpha$

\subsection{Lions' Lemma}
Let $\mathcal O$ be an open, bounded and Lipschitz-regular subset of $\mathbb R^{d+1}\approx\R_t\times\R^d_x$.
We recall that 
$$\mathcal{H}^{-1}(\mathcal O) = \left\{ w \in \mathcal D^\ast(\mathcal O)\,:\,|\langle w,u\rangle_{\!_{\mathcal O}}| \leq C\,\|u\|_{\mathcal{H}^1_0(\mathcal O)},\,C>0 \right\},$$
where $\mathcal{D}^\ast(\mathcal O)$ denotes the space of distributions over $\mathcal O,$ equipped with the weak$ \ast$ topology, and $\langle w,u\rangle_{\!_{\mathcal O}}$ is the duality product between $\mathcal{H}^{-1}$ and $\mathcal{H}^1_0.$
The norm on~$\mathcal{H}^1_0(\mathcal O)$ is as usual $ u \mapsto \|\nabla u\|_{\mathrm{L}^2(\mathcal O)}.$ 
On $\mathcal{H}^1(\mathcal O),$ we introduce the norm
$$\|u\|^2_{\mathrm{H}^1} = \left| \int_{\mathcal O} u\,dt\,dx \right|^2 + \|\nabla u\|^2_{\mathrm{L}^2(\mathcal O)}.$$
The norm induced on $\mathrm{H}^{-1}(\mathcal O)$ is then
$$\|w\|_{\mathrm{H}^{-1}(\mathcal O)}^2 = \langle w,1\rangle_{\!_{\mathcal O}}^2 + \|z_w\|^2_{\mathrm{H}^1(\mathcal O)},$$
where $z_w$ is the solution to 
$$ -(\partial_{tt} + \Delta_x)\,z_w = w - \langle w,1\rangle_{\!_{\mathcal O}}, \quad \int_{\mathcal O} z_w\,dt\,dx = 0.$$
Lions' Lemma gives a sufficient condition for a distribution to be an $\mathrm{L}^2$ function. The following statement is taken from~\cite{amrouche2015lemma}.
\begin{lemma}\label{lemma2}
Let $\mathcal O$ be a bounded, open and Lipschitz-regular subset in $\mathbb R^{d+1}.$ Then, for all $u \in \mathcal{D}^\ast(\mathcal O),$ we have that $u \in \mathrm{L}^2(\mathcal O)$ if and only if the weak gradient $\nabla u$ belongs to $\mathrm{H}^{-1}(\mathcal O).$ Moreover, there exists a constant $C_L(\mathcal O)$ such that
\begin{equation*}
 \left\|u - \int_{\mathcal O} u\,dx\,dt \right\|^2_{\mathrm{L}^2(\mathcal O)} \leq C_L \|\nabla u\|^2_{H^{-1}(\mathcal O)},
\end{equation*}
for any $u \in \mathrm{L^2}(\mathcal O).$
\end{lemma}
According to~\cite{calderon1956singular,bogovskii1979solution,csato2011pullback}, if $\mathcal O$ is star-shaped w.r.t.~a ball, then the constant~$C_L$ has the following structure:
\begin{equation}
 \label{eq20}
 C_L = 4\,|S^d|\,\frac{\mathrm{D}(\mathcal O)}{\mathrm d(\mathcal O)},
\end{equation}
where $\mathrm{D}$ is the diameter of $\mathcal O$, while $\mathrm d(\mathcal O)$ is the diameter of the largest ball one can include in $\mathcal O.$ See in particular ~\cite[Remark 9.3]{csato2011pullback} and~\cite[Lemma 1]{bogovskii1979solution}. As a consequence, we have the following explicit expression of~$C_L$ when $\mathcal O=\Omega.$
\begin{lemma}\label{bestlions}
Let $L>0$, $\tau \in (0,L)$ and $\Omega = (0,\tau) \times (0,L)^d.$ Lemma~\ref{lemma2} holds with
\begin{equation}\label{bestlionseq}
C_L = 4\,|S^d|\,\frac{\sqrt{d\,L^2+\tau^2}}{\tau}. 
\end{equation}
\end{lemma}

\subsection{The kinetic Ornstein-Uhlenbeck equation}

We consider solutions to~\eqref{eq2} in the weak sense, \emph{i.e.}, functions $ h$ in the space $ \mathrm C(\R^+; \mathrm{L}^2(dx\,d\gamma_\alpha))$ with initial datum $h_0 = h(0,\cdot,\cdot)$ in $\mathrm{L}^2(dx\,d\gamma_\alpha)$ such that~\eqref{eq2} holds in the sense of distributions on $(0,\infty)\times\R^d_x\times\R^d_v.$ 
The following result is taken from~\cite{albritton2019variational} if $\alpha=2$. The extension to $\alpha\neq2$ is straightforward as follows from a careful reading of the proof in~\cite[Proposition~6.10]{albritton2019variational}.
\begin{proposition}\label{Prop:Existence}
Let $L>0$ and $\alpha >0.$ With $\Omega = (0,\tau) \times (0,L)^d,$ for all zero-average initial datum $h_0 \in \mathrm{L}^2(dx\,d\gamma_\alpha)$, there exists a unique solution $h$ to~\eqref{eq2} such that $h \in \mathrm H_\kin$ for all $\tau > 0.$
\end{proposition}
Regularity properties for~\eqref{eq2} are collected in~\cite[Section 6]{albritton2019variational}. In the special case \hbox{$\alpha = 2$,} some fractional regularity along all directions of the phase space are known. Also see~\cite{perthame2004mathematical} for further result on regularity theory for kinetic Fokker-Planck equations.

\subsection{A priori estimates}
We state two estimates for solutions to~\eqref{eq2}.
\begin{lemma}\label{lemma3} Let $L>0$, $\tau >0,$ $\Omega = (0,\tau) \times (0,L)^d,$ and $\alpha >0.$ If $h$ is a solution to~\eqref{eq2}, then we have
\begin{equation}\label{eq21}
\|(\partial_t + v \cdot \nabla_x)\,h \|_{\mathrm L^2(\Omega;{\mathrm{H}^{-1}_\alpha})} \leq \| \nabla_v h \|_{\mathrm{L^2}(\Omega;\mathrm{L}^2_\alpha)}.
\end{equation}
\end{lemma}
\begin{proof}
Take a test function $\phi \in \mathrm{L}^2(\mathrm{H}^1_\alpha),$ and write 
$$\int_{\mathbb T} \langle (\partial_t + v \cdot \nabla_x)\,h , \phi \rangle\,dx = \int_{\mathbb T} \langle \Delta_\alpha h , \phi \rangle\,dx = -\int_{\mathbb T}(\nabla_v h, \nabla_v \phi)\,dx,$$
from which~\eqref{eq21} easily follows, after maximizing over $\|\nabla_v\phi\|_{\mathrm L^2_\alpha}\le1.$
\end{proof}
For completeness, let us recall the classical $\mathrm L^2$ decay estimate for solutions to~\eqref{eq2}.
\begin{lemma}\label{energylem}
Let $L>0$, $\tau >0,$ $\Omega = (0,\tau) \times (0,L)^d,$ and $\alpha >0.$ If $h$ is a solution to~\eqref{eq2}, then we have
\begin{equation*}
 \frac{d}{dt} \|h\|^2_{\mathrm{L}^2(dx\,d\gamma_\alpha)} = -\,2\,\|\nabla_vh\|^2_{\mathrm{L}^2(dx\,d\gamma_\alpha)}.
\end{equation*}
\end{lemma}

\section{An averaging lemma and a generalized Poincar\'e inequality}\label{Sec3}
For all functions $h \in \mathrm{H}_\kin,$ we define the \emph{spatial density}
$$\rho_h := \int_{\R^d} h(\cdot,\cdot,v)\,d\gamma_\alpha.$$
Notice that $\int_Q\rho_h\,dx=0$ whenever $h$ is a zero-average function.

\subsection{Averaging lemma}
Inspired by~\cite[Proposition 6.2]{albritton2019variational}, the following \emph{averaging lemma} provides a norm of the spatial density, as for instance in~\cite{perthame2004mathematical}.
\begin{lemma}\label{lemmaavg}
Let $L>\tau >0,$ $\Omega = (0,\tau) \times (0,L)^d,$ and $\alpha >0.$ For all $h \in \mathrm{H}_\kin,$ we have
\begin{equation}\label{eq7}
\|\nabla_{t,x}\rho_h\|^2_{\mathrm{H}^{-1}(\Omega)} \leq d_\alpha \left( \|h-\rho_h\|^2_{\mathrm{L^2}(dt\,dx\,d\gamma_\alpha)} + \|\partial_t h +v\cdot\nabla_x h\|^2_{\mathrm{L}^2(\Omega; H^{-1}_\alpha)} \right)
\end{equation}
with
\begin{equation}\label{dalpha}
d_\alpha = 2\,\left(\|v_1 |v|^2\|^2_{\mathrm{L}^2_\alpha} + \left(1+ \tfrac{L^2}{4\,\pi^2}\right)\||v|^2\|^2_{\mathrm{L}^2_\alpha} + \tfrac{d^2L^2}{4\,\pi^2}\,\|v\|^2_{\mathrm{L}^2_\alpha} \right).
\end{equation}
\end{lemma}
Inequality~\eqref{eq7} can be extended to any measure $d\gamma$ such that $\int_{\R^d}|v|^4\,d\gamma<\infty$ and $\int_{\R^d}v\,d\gamma=0$.
The proof of Lemma~\ref{lemmaavg} is technical, but follows in a standard way from the time-independent case, as it is common in averaging lemmas: see ~\cite{perthame2004mathematical}.
For sake of simplicity, we detail only the $t$-independent case below.

\begin{proof}[Proof of Lemma~\ref{lemmaavg}]
Assume that $h \in \mathrm H_\kin$ does not depend on $t.$ 
Let $\phi \in \mathcal{D}(\mathbb T)^{d}$ be a smooth test-vector field with zero average on each component. 
We write
$$- \int_{\mathbb T} \rho_h \nabla_x \cdot \phi\,dx = \int_{\mathbb T} (\nabla_x \rho_h) \cdot \phi\,dx ,$$
with a slight abuse of notation, since the integral of the r.h.s.~is in fact a duality product. 
Using $\int_{\mathbb R^d} v_i\,v_j\,d\gamma_\alpha = d^{-1}\,\|v\|_{\mathrm{L}^2_\alpha}^2\,\delta_{ij},$ we obtain
$$\int_{\mathbb T} \nabla_x \rho_h \cdot \phi\,dx= d\,\|v\|_{\mathrm{L}^2_\alpha}^{-2}\iint_{\mathbb T \times \mathbb R^d} v \cdot \nabla_x \rho_h\,\phi \cdot v\,dx\,d\gamma_\alpha.$$
By adding and subtracting $\rho_h$, and then integrating by parts, still at formal level, we obtain
\begin{multline*}
\iint_{\mathbb T \times \mathbb R^d} v \cdot \nabla_x \rho_h\,\phi \cdot v\,dx\,d\gamma_\alpha\\
=-\iint_{\mathbb T \times \mathbb{R}^d} v \cdot (h - \rho_h)\,\nabla_x\phi \cdot v\,dx\,d\gamma - \iint_{\mathbb T \times \mathbb{R}^d} v \cdot \nabla_x h\,v \cdot \phi\,\,dx\,d\gamma_\alpha\\
\le\| h - \rho_h \|_{\mathrm{L}^2(dx\,d\gamma_\alpha)}\,\|\nabla_x\phi\|_{\mathrm{L}^2(dx)} \left\|v^2\right\|_{\mathrm{L}^2_\alpha}\\
+\|v \cdot \nabla_x h\|_{\mathrm{L}^2(\mathrm{H}^{-1}_\alpha)}\,\|\phi\|_{\mathrm{L}^2(dx)}\,\|v\|_{\mathrm{L}^2_\alpha}
\end{multline*}
using Cauchy-Schwarz inequalities and duality estimates. By the Poincar\'e inequality, we know that
$$\frac{4\,\pi^2}{d\,L^2}\,\|\phi\|_{\mathrm{L}^2(dx)}^2\le\|\nabla\phi\|_{\mathrm{L}^2(dx)}^2
$$
Maximizing the r.h.s. on~$\phi$ such that $\|\nabla_x\phi\|_{\mathrm{L}^2(dx)}\le1$ completes the proof of the $t$-independent case. When $h$ additionally depends on $t$, the same scheme can be applied with $v\cdot\nabla_x$ replaced by  $\partial_t+v\cdot\nabla_x$.
\end{proof}

\subsection{A generalized Poincar\'e inequality}
The next \emph{a priori} estimate is at the core of the method. It is a modified Poincar\'e inequality in $t$, $x$ and $v$ which relies on Lemma~\ref{lemmaavg} and involves derivatives of various orders.
\begin{proposition}\label{Pro:pgPI}
Let $L>0$, $\tau >0,$ $\Omega = (0,\tau) \times (0,L)^d,$ and $\alpha >0.$ 
Then, for all $h \in \mathrm{H}_\kin$ with zero average, we have that 
\begin{equation}
 \label{poigen}
 \|h\|^2_{\mathrm{L^2}(dt\,dx\,d\gamma_\alpha)} \leq C \left( \|h-\rho_h\|^2_{\mathrm{L^2}(dt\,dx\,d\gamma_\alpha)} + \|\partial_t h +v\cdot\nabla_x h\|^2_{\mathrm{L}^2(\Omega; H^{-1}_\alpha)} \right)
\end{equation}
with $C = 1+C_L\,d_\alpha$, where $C_L$ and $d_\alpha$ are given respectively by~\eqref{bestlionseq} and~\eqref{dalpha}.
\end{proposition}
\begin{proof}
By orthogonality in $\mathrm{L}^2(\Omega;\mathrm{L}^2_\alpha)$ and because $d\gamma_\alpha$ is a probability measure, we have the decomposition
$$\|h\|^2_{\mathrm{L}^2(\Omega;\mathrm{L}^2_\alpha)} = \|h-\rho_h\|^2_{\mathrm{L}^2(\Omega;\mathrm{L}^2_\alpha)} + \|\rho_h\|^2_{\mathrm{L}^2(\Omega)}.$$
The function $\rho_h$ has zero average on $\Omega$ by construction, so that 
$$\|h\|^2_{\mathrm{L}^2(\Omega;\mathrm{L}^2_\alpha)} \leq \|h-\rho_h\|^2_{\mathrm{L}^2(\Omega;\mathrm{L}^2_\alpha)} + C_L\,\|\nabla_{x,t} \rho\|^2_{\mathcal{H}^{-1}(\Omega)}$$
by Lemma~\ref{bestlions}. Hence
$$\|h\|^2_{\mathrm{L}^2(\Omega;\mathrm{L}^2_\alpha)} \leq (1 + C_L\,d_\alpha)\,\|h-\rho_h\|^2_{\mathrm{L^2}(\Omega;d\gamma_\alpha)}+ C_L\,d_\alpha\,\|\partial_t h +v\cdot\nabla_x h\|^2_{\mathrm{L}^2(\Omega; H^{-1}_\alpha)}$$
by Lemma~\ref{lemmaavg}. This concludes the proof.
\end{proof}

\section{Linear and superlinear local equilibria: exponential decay rate}\label{Sec4}

In this Section, we consider the case $\alpha\ge1$ and the domain $\Omega = (t,t+\tau) \times (0,L)^d,$ for an arbitrary $t\ge0.$ Let us define $\kappa_\alpha := (1+C_L\,d_\alpha)(P_\alpha+1)$
where $P_\alpha$ is the Poincar\'e constant in~\eqref{poi} and where $C_L$ and $d_\alpha$ are given respectively by~\eqref{bestlionseq} and~\eqref{dalpha}.
\begin{lemma}\label{poigeneral}
Let $L>0$, $\tau >0,$ $t\ge0,$ $\Omega_t = (t,t+\tau) \times (0,L)^d,$ and $\alpha\ge1.$ Then, for all $h \in \mathrm{H}_\kin$ with zero average which solve \eqref{eq2}, we have that
\begin{equation}\label{genpoinc2}
 \|h\|^2_{\mathrm{L^2}(dt\,dx\,d\gamma_\alpha)} \leq \kappa_\alpha\,\|\nabla_v h\|^2_{\mathrm{L^2}(dt\,dx\,d\gamma_\alpha)}.
 \end{equation}
\end{lemma}
\begin{proof} We know that
\begin{equation*}
 \|h\|^2_{\mathrm{L^2}(dt\,dx\,d\gamma_\alpha)} \leq (1+C_L\,d_\alpha) \left(P_\alpha\,\|\nabla_v h\|^2_{\mathrm{L^2}(dt\,dx\,d\gamma_\alpha))} + \|\partial_t h +v\cdot\nabla_x h\|^2_{\mathrm{L}^2(H^{-1}_\alpha)} \right)
 \end{equation*}
as a consequence of~\eqref{poi} and~\eqref{poigen}. Then~\eqref{genpoinc2} follows from Lemma~\ref{lemma3}.
\end{proof}

We are ready to prove Theorem~\ref{thm1} with an explicit estimate of the constant~$\lambda.$
\begin{proof}[Proof of Theorem \ref{thm1}]
Inequality \eqref{genpoinc2} -- on the interval $(t,t+\tau)$ -- gives
$$\int_t^{t+\tau }\|h(s,\cdot,\cdot)\|^2_{\mathrm{L^2}(dx\,d\gamma_\alpha)}\,ds \leq \kappa_\alpha \int_t^{t+\tau} \|\nabla_v h(s,\cdot,\cdot)\|^2_{\mathrm{L^2}(dx\,d\gamma_\alpha)}\,ds.$$
With $\lambda=2/\kappa_\alpha,$ we deduce from Lemma~\ref{energylem} that
\begin{multline*}
\frac{d}{dt} \int_t^{t+\tau} \| h(s,\cdot,\cdot)\|^2_{\mathrm{L^2}(dx\,d\gamma_\alpha)}\,ds = -\,2\int_t^{t+\tau}\|\nabla_v h(s,\cdot,\cdot)\|^2_{\mathrm{L}^2(dx\,d\gamma_\alpha)}\,ds\\
\le-\,\lambda\int_t^{t+\tau} \| h(s,\cdot,\cdot)\|^2_{\mathrm{L^2}(dx\,d\gamma_\alpha)}\,ds.
\end{multline*}
Gr\"onwall's Lemma
and the monotonicity of $t\mapsto\| h(t,\cdot,\cdot)\|^2_{\mathrm{L^2}(dx\,d\gamma_\alpha)}$ imply
\begin{multline*}\int_t^{t+\tau }\|h(s,\cdot,\cdot)\|^2_{\mathrm{L^2}(dx\,d\gamma_\alpha)}\,ds \leq\int_0^\tau\|h(s,\cdot,\cdot)\|^2_{\mathrm{L^2}(dx\,d\gamma_\alpha)}\,ds\,e^{-\lambda\,t}\\
\leq \tau\,\|h_0\|^2_{\mathrm{L}^2(dx\,d\gamma_\alpha)}\,e^{-\lambda\,t}
\end{multline*}
for ant $t\ge0$, which proves~\eqref{eq5}, that is, Theorem~\ref{thm1}.
\end{proof}
Indeed, the estimate for $\lambda$ is explicit, as we state in the following. 
\begin{proposition}\label{thm1vero}
For any $\alpha \geq 1$, Theorem~\ref{thm1} holds true with 
$$\frac1\lambda = \frac1\tau\left(\tau + \sqrt{d\,L^2 + \tau^2}\right)\Big(2\,d_\alpha\,|\mathbb S^{d-1}|\,(P_\alpha +1)\Big).$$
\end{proposition}
Notice that the r.h.s.~vanishes as $\tau \downarrow 0,$ which is expected because of the degeneracy of $\Delta_\alpha$: an exponential decay rate of $\| h(t,\cdot,\cdot)\|^2_{\mathrm{L^2}(dx\,d\gamma_\alpha)}$ cannot hold. 
The section is concluded showing how the result above yields the classic hypocoercivity estimate of Corollary \ref{cor1}.
\begin{proof}[Proof of Corollary~\ref{cor1}]
For any $t\ge0$, we know from Theorem~\ref{thm1} that
$$\|h(t+\tau,\cdot,\cdot) \|^2_{\mathrm{L}^2(dx\,d\gamma_\alpha)}\le\strokedint_t^{t+\tau} \|h(s,\cdot,\cdot) \|^2_{\mathrm{L}^2(dx\,d\gamma_\alpha)}\,ds \leq\,\|h_0\|^2_{\mathrm{L}^2(dx\,d\gamma_\alpha)}\,e^{-\lambda\,t},$$
as a consequence of the monotonicity of the $\mathrm{L}^2$ norm, according to Lemma~\ref{energylem}. This proves that
$$\|h(t,\cdot,\cdot) \|^2_{\mathrm{L}^2(dx\,d\gamma_\alpha)}\le\,C\,\|h_0\|^2_{\mathrm{L}^2(dx\,d\gamma_\alpha)}\,e^{-\lambda\,t}$$
with $C=e^{\lambda\,\tau}$ for any $t\ge\tau$. However, if $t\in[0,\tau)$, it turns out that $C\,e^{-\lambda\,t}\ge 1$ so that the inequality is also true by Lemma~\ref{energylem}.
This concludes the proof.
\end{proof}

\section{Sublinear equilibria: algebraic decay rates}\label{Sec5}

\subsection{Proof of Theorem~\texorpdfstring{\ref{thm2}}{3}}
Assume that $\alpha\in (0,1).$ Let us define the parameter $\beta = 2\,(1-\alpha)/p$ where $p$, $q>1$ are H\"older conjugate exponents, \emph{i.e.}, $\frac1p+\frac1q=1$ and define
\begin{equation}\label{Zh}
Z_h(t) :=\iint_{\mathbb T\times\mathbb{R}^d}\langle v\rangle^{\beta\,q}\,|h-\rho_h|^2\,dx\,d\gamma_\alpha.
\end{equation}
The following estimates replace Proposition~\ref{Pro:pgPI}. 
\begin{proposition}\label{prop3}
Let $L>0$, $\tau >0,$ $t\ge0,$ $\Omega_t = (t,t+\tau) \times (0,L)^d,$ and $\alpha\ge1.$ With the above notations, for  all $h \in \mathrm{H}_\kin$ with zero average, we have that 
\begin{multline*}
 \|h\|^2_{\mathrm{L^2}(dt\,dx\,d\gamma_\alpha)}\\ \leq C\,P_\alpha^\frac1p\,\|\nabla_v h\|^\frac2p_{\mathrm{L}^2(dt\,dx\,d\gamma_\alpha)} \left({\textstyle\int_t^{t+\tau}Z_h(s)\,ds}\right)^{\frac{1}{q}} + C\,\|\partial_t h +v\cdot\nabla_x h\|^2_{\mathrm{L}^2(\Omega; H^{-1}_\alpha)}
\end{multline*}
where $C = 1+C_L\,d_\alpha$ is as in Proposition~\ref{Pro:pgPI} and $P_\alpha$ denotes the constant in the weighted Poincar\'e inequality~\eqref{weightpeq}.\end{proposition}
\begin{proof} Using~\eqref{poigen} and H\"older's inequality w.r.t.~the variable $v$, we find that
$$\|h-\rho_h\|_{\mathrm L^2_\alpha}^2\le\left(\int_{\R^d}\langle v\rangle^{-\beta\,p}\,|h-\rho_h|^2\,d\gamma_\alpha \right)^{\frac{1}{p}} \left( \int_{\R^d}\langle v\rangle^{\beta\,q}\,|h-\rho_h|^2\,d\gamma_\alpha \right)^{\frac{1}{q}}.$$
The weighted Poincar\'e inequality~\eqref{weightpeq} with $\beta\,p=2\,(1-\alpha)$ and an additional H\"older inequality w.r.t.~the variables $t$ and $x$ allow us to complete the proof.\end{proof}
\begin{lemma}\label{lem4}
Let $L>0$, $\tau >0,$ $t\ge0,$ $\Omega_t = (t,t+\tau) \times (0,L)^d,$ and $\alpha\in(0,1).$ There is a constant $W>0$ such that, for all solution $h \in \mathrm{H}_\kin$ to~\eqref{eq2} with an initial datum $h_0$ with zero average,  using the notation~\eqref{Zh} as in Proposition~\ref{prop3}, we have
$$Z_h(t)\le W\iint_{\mathbb T\times\R^d}\langle v\rangle^{\beta\,q}\,h_0^2\,dx\,d\gamma_\alpha,\quad\forall\,t\ge0.$$
\end{lemma}
\begin{proof}
An elementary computation shows that
\begin{align*}
\iint_{\mathbb T\times\R^d}\langle v\rangle^{\beta\,q}\,|h-\rho_h|^2\,dx\,d\gamma_\alpha&\le2\iint_{\mathbb T\times\R^d}\langle v\rangle^{\beta\,q}\left(h^2+\rho_h^2\right)\,dx\,d\gamma_\alpha\\
&\le2\left({\textstyle 1+\int_{\R^d}\langle v\rangle^{\beta\,q}\,d\gamma_\alpha}\right)\iint_{\mathbb T\times\R^d}\langle v\rangle^{\beta\,q}\,h^2\,dx\,d\gamma_\alpha
\end{align*}
because $\rho_h^2=\left(\int_{\R^d}h\,d\gamma_\alpha\right)^2\le\int_{\R^d}h^2\,d\gamma_\alpha\le\int_{\R^d}\langle v\rangle^{\beta\,q}\,h^2\,d\gamma_\alpha.$ According to~\cite[Proposition~4]{bouin2019hypocoercivity}, there is a constant $\mathcal K_{\beta\,q}>1$ such that
$$\iint_{\mathbb T\times\R^d}\langle v\rangle^{\beta\,q}\,|h(t,x,v)|^2\,dx\,d\gamma_\alpha\le\mathcal K_{\beta\,q}\,\iint_{\mathbb T\times\R^d}\langle v\rangle^{\beta\,q}\,h_0^2\,dx\,d\gamma_\alpha,\quad\forall\,t\ge0.$$
The result follows with $W=2\left( 1+\int_{\R^d}\langle v\rangle^{\beta\,q}\,d\gamma_\alpha\right)\mathcal K_{\beta\,q}.$
\end{proof}

Assume that $h \in \mathrm{H}_\kin$ solves~\eqref{eq2} with an initial datum $h_0$ with zero average and let us collect our estimates. With Proposition~\ref{lemma3}, Proposition~\ref{prop3}, and Lemma~\ref{lem4}, the estimate of Lemma~\ref{poigeneral} is replaced by
\begin{equation}\label{genpoinc2b}\|h\|^2_{\mathrm{L^2}(dt\,dx\,d\gamma_\alpha)}\le A\,\|\nabla_v h\|^\frac2p_{\mathrm{L}^2(dt\,dx\,d\gamma_\alpha)} + C\,\|\nabla_v h\|^2_{\mathrm{L}^2(dt\,dx\,d\gamma_\alpha)}\end{equation}
with $A=C\,P_\alpha^\frac1p \left(\tau\,W\right)^{1/q}\left(\iint_{\mathbb T\times\R^d}\langle v\rangle^{\beta\,q}\,h_0^2\,dx\,d\gamma_\alpha\right)^{1/q}.$

\medskip The main result of the section is a technical version of Theorem~\ref{thm2}. Let
$$\mathsf x(t):=\strokedint_t^{t+\tau} \|h(s,\cdot,\cdot)\|^2_{\mathrm{L}^2(dx\,d\gamma_\alpha)}\,ds\quad\mbox{and}\quad\mathsf y(t):=\strokedint_t^{t+\tau} \|\nabla_vh(s,\cdot,\cdot)\|^2_{\mathrm{L}^2(dx\,d\gamma_\alpha)}\,ds$$
where norms are taken on $\mathbb T\times\R^d$. We know from Lemma~\ref{energylem} and~\eqref{genpoinc2b} that
$$\mathsf x'=-\,2\,\mathsf y\quad\mbox{and}\quad\mathsf x\le\varphi(\mathsf y):=A\,\mathsf y^{1/p}+C\,\mathsf y.$$
Finally, let us denote by $\varphi^{-1}$ the inverse of $\mathsf y\mapsto\varphi(\mathsf y)$ and consider 
$$\psi(\mathsf z):=\int_{\mathsf z}^{\mathsf x_0}\frac{dz}{2\,\varphi^{-1}(z)}\quad\mbox{with}\quad\mathsf x_0=\|h_0\|^2_{\mathrm{L^2}(dx\,d\gamma_\alpha)}.$$
\begin{theorem}\label{thm2vero}
Let $L>0$, $\tau >0,$ $t\ge0,$ $\Omega_t = (t,t+\tau) \times (0,L)^d,$ and $\alpha\in(0,1).$ With the above notations, for all solution $h \in \mathrm{H}_\kin$ to~\eqref{eq2} with an initial datum~$h_0$ with zero average, we have
$$\mathsf x(t)=\strokedint_t^{t+\tau} \|h(s,\cdot,\cdot)\|^2_{\mathrm{L}^2(dx\,d\gamma_\alpha)}\,ds\le\psi^{-1}(t),\quad\forall\,t\ge0.$$
\end{theorem}
\begin{proof} The strategy goes as in~\cite{liggett1991l2,bouin2019hypocoercivity}. Everything reduces to the differential inequality
$$\mathsf x'\le-\,2\,\varphi^{-1}(\mathsf x)$$
using the monotonicity of $\mathsf y\mapsto\varphi(\mathsf y)$. From by the elementary Bihari-Lasalle inequality, see~\cite{bihari1956generalization,lasalle1949uniqueness}, which is obtained by a simple integration, we obtain
$$\mathsf x(t)=\strokedint_t^{t+\tau} \|h(s,\cdot,\cdot)\|^2_{\mathrm{L}^2(dx\,d\gamma_\alpha)}\,ds\le\psi^{-1}\big(t+\psi(\mathsf x(0))\big).$$
Since, on the one hand
$$\mathsf x(0)=\strokedint_0^\tau\|h(s,\cdot,\cdot)\|^2_{\mathrm{L}^2(dx\,d\gamma_\alpha)}\,ds\le\|h_0\|^2_{\mathrm{L}^2(dx\,d\gamma_\alpha)}=\mathsf x_0$$
because $s\mapsto\|h(s,\cdot,\cdot)\|^2_{\mathrm{L}^2(dx\,d\gamma_\alpha)}$ is nonincreasing according to Lemma~\ref{energylem}, and $\psi$ is nonincreasing on the other hand, then
$$\psi^{-1}\big(t+\psi(\mathsf x(0))\big)\le\psi^{-1}(t),$$
which concludes the proof. Notice that the dependence on $h_0$ enters in $A$ and~$\mathsf x_0$, and henceforth in $\varphi$ and $\psi$.\end{proof}

\begin{proof}[Proof of Theorem~\ref{thm2}] Since $\lim_{t\to+\infty}\mathsf y(t)=0$, we have that $\varphi(\mathsf y(t))\sim A\,\mathsf y(t)^{1/p}$ as $t\to+\infty$, which heuristically explains the role played by $p$ in~\eqref{decay}. This can be made rigorous as follows. Notice that
$$\varphi(\mathsf y)=A\,\mathsf y^{1/p}+C\,\mathsf y\le A_0\,\mathsf y^{1/p},\quad\forall\,\mathsf y\le\mathsf y_0,\quad\mbox{with}\quad A_0=A+C\,\mathsf y_0^{1-1/p}.$$
With $A$ replaced by $A_0$ and $C$ replaced by $0$, the computation of the proof of Theorem~\ref{thm2vero} is now explicit. With the choice $\mathsf y_0=\varphi^{-1}(\mathsf x_0),$ we know that $\mathsf y(t)\le\mathsf y_0$ for any $t\ge0$ and obtain
\be{alpha01}\mathsf x(t)\le\left(\mathsf x_0^{1-p}+2\,(p-1)\,A_0^{-p}\,t\right)^{-\frac1{p-1}},\quad\forall\,t\ge0.\ee
Using $\mathsf x_0=A\,\mathsf y_0^{1/p}+C\,\mathsf y_0\ge C\,\mathsf y_0$, we know that $A_0=\mathsf x_0\,\mathsf y_0^{-1/p}\le A+C^{1/p}\,\mathsf x_0^{1-1/p}$, which proves~\eqref{decay} with
$$K=\max\left\{1,\big(2\,(p-1)\big)^{1/(1-p)}\left(C\,P_\alpha^{1/p} \left(\tau\,W\right)^{1-1/p}+C^{1/p}\right)\right\}.$$
The conclusion holds using $\sigma=\beta\,q=2\,(1-\alpha)/(p-1).$
\end{proof}

\subsection{The linear threshold: from algebraic to exponential rates}
A very natural question arises: \emph{is the result Theorem~\ref{thm1vero} (corresponding to $\alpha\in(0,1)$ consistent with the result of Theorem~\ref{thm2vero} (which covers any $\alpha\ge1$) ?} A first observation is that we can vary $\alpha$ in the assumptions concerning the initial data.
\begin{lemma} If $h_0\in\mathrm L^2(\mathbb T;\mathrm L^2_{\alpha_0})$ for some $\alpha_0\in(0,1)$, then $\langle v\rangle^{\sigma/2}\,h_0\in\mathrm L^2(\mathbb T;\mathrm L^2_\alpha)$ for any $\alpha>\alpha_0$ and any $\sigma>0.$
\end{lemma}
The proof is a simple consequence of the fact that $v\mapsto\langle v\rangle^\sigma\,\exp\langle v\rangle^{\alpha-\alpha_0}$ is uniformly bounded. For any $\alpha>\alpha_0$, let us denote the corresponding solution of~\eqref{eq2} with initial datum $h_0$, of zero average, by $h^{(\alpha)}$.

If $\alpha\in(\alpha_0,1)$, then~\eqref{alpha01} can be rewritten as
$$\strokedint_t^{t+\tau}\|h^{(\alpha)}(s,\cdot,\cdot)\|^2_{\mathrm{L}^2(dx\,d\gamma_\alpha)}\,ds\le\|h_0\|^2_{\mathrm{L^2}(dx\,d\gamma_\alpha)}\,\big(1+(p-1)\,\ell(\alpha)\,t\big)^{-\frac1{p-1}}.
$$
By passing to the limit as $\alpha\to1^-$, we recover~\eqref{eq5} with $\lambda=\lim_{\alpha\to1^-}\ell(\alpha)$, where $\ell(\alpha)=2\,\|h_0\|^{2(p-1)}_{\mathrm{L^2}(dx\,d\gamma_\alpha)}\,A_0^{-p}$ and $A_0=A_0(\alpha)$ as above. The Poincar\'e constant~$P_\alpha$ in the weighted Poincar\'e inequality~\eqref{weightpeq} admits a limit as $\alpha\to1^-$, according to~\cite[Appendix~A]{bouin2017hypocoercivity}. 

The limit of~$\lim_{\alpha\to1^-}\ell(\alpha)$ is certainly not optimal. By working directly on the Bihari-Lasalle estimate of Theorem~\ref{thm2vero}, we can recover the value of $\lambda$ in Theorem~\ref{thm1vero}. Notice here that $\sigma>0$ plays essentially no role and can be taken arbitrarily small, even depending on $\alpha$, but such that $p=2\,(1-\alpha)/\sigma\to1$ as $\alpha\to1^-.$

As in \cite{MR3625186,MR4265692}, it is possible to obtain improved decay rates in \eqref{decay} by picking the initial datum in a smaller space. Typically, the control of additional norms or moments is asked. However, the strategy in the current paper is in the opposite direction. If $\alpha \in (0,1)$ we are interested in taking the initial data in a space as large as possible so that we can compute decay rates. The additional conditions to be imposed have been shown to vanish as $\alpha \to 1^-.$
\section{Hypocoercivity and comparison with some other methods}\label{Sec7}

\subsection{An explicit hypocoercivity result}
Theorem~\ref{thm1} implies an $\mathrm L^2-$hypocoercivity result in the linear and superlinear regimes $\alpha \geq 1,$ see Corollary \ref{cor1}. 
The remainder of this section is devoted to a comparison with earlier  hypocoercivity results in a simple benchmark case: let $\alpha=2$, $d = 1$ and $L = 2\pi.$ In this case, for the choice $\tau = 2\pi,$ Theorem~\ref{thm1} amounts to
$$\|h\|^2_{\mathrm{L}^2(dx\,d\gamma_2)} \leq \mathrm{e}^{\frac\pi{8\,\sqrt3}}\,\|h_0\|^2_{\mathrm{L}^2(dx\,d\gamma_2)}\,\mathrm{e}^{-\frac{t}{8\,\sqrt3}},\quad\forall\,t\ge0.$$
For sake of comparison, notice that $\lambda=1/(8\,\sqrt3)\approx0.0721688$.
Even if we are aware of explicit or sharp results in other metrics than $\mathrm{L}^2$ for \eqref{eq2}, as \cite{dietert2015convergence,Mouhot_2006}, we restrict our discussion to $\mathrm{L}^2-$hypocoercivity methods.
\subsection{The DMS method}
The first comparison is with the abstract twisted $\mathrm L^2$ hypocoercivity method of~\cite{dolbeault2015hypocoercivity,bouin2017hypocoercivity}. Let $\|\cdot\|$ be the norm of $\mathrm{L}^2(dx\,d\gamma_2)$ and $(\cdot,\cdot)$ the associated scalar product. We consider the evolution equation
\be{Eqn:abstract}\partial_t h + Th = \mathcal L\,h.\ee
\begin{theorem}\label{Thm:DMS} Let $h$ be a solution of~\eqref{Eqn:abstract} with initial datum $h_0\in\mathrm L^2(Q;\mathrm L^2_2)$ and assume that $T$ and $\mathcal L$ are respectively anti-self-ajoint and self-adjoint operators on $\mathrm L^2(Q;\mathrm L^2_2)$ such that, for some positive constants $\lambda_m$, $\lambda_M$ and $C_M$, we have
\begin{enumerate}
\item[(A1)] $(-\mathcal Lh,h) \geq \lambda_m\,\|(1-\Pi)\,h\|^2$ for all $h \in D(\mathcal L),$
\item[(A2)] $\|T\,\Pi\,h \|^2\geq \lambda_M\,\|\Pi\,h\|^2$ for all $h \in D(T\,\Pi),$
\item[(A3)] $\Pi\,T\,\Pi\,h = 0,$
\item[(A4)] $\|A\,T\,(1-\Pi)\,h\| + \|A\,\mathcal L\,h\| \leq C_M\,\|(\mathrm{Id}-\Pi)\,h\|$ for all admissible $h\in\mathrm L^2(Q;\mathrm L^2_2)$
\end{enumerate}
where $A:=\big(\mathrm{Id}+(T\,\Pi)^\ast\,T\,\Pi\big)^{-1}\,(T\,\Pi)^\ast$ and $\Pi$ is the projection in $\mathrm L^2_2$ onto the kernel of $\mathcal L$. Then we have
$$\|h(t,\cdot,\cdot)\|^2\le C\,\|h_0\|^2\,e^{-\lambda\,t}\quad\forall\,t\ge0$$
with $C=(1+\delta)/(1-\delta),$ $\delta=\frac12\,\min\big\{1,\lambda_m,\tfrac{\lambda_m\,\lambda_M}{(1+\lambda_M)\,C_M^2}\big\}$ and $\lambda=\frac{2\,\delta\,\lambda_M}{3\,(1+\lambda_M)}.$
\end{theorem}
This result is taken from~\cite[Proposition~4]{bouin2017hypocoercivity}. According to~\cite[Corollary~9]{bouin2017hypocoercivity}, we have the estimate $\lambda = 1/24\approx0.041667.$ A minor improvement is obtained as follows. Using  Theorem~\ref{Thm:DMS} applied with $d=1$, $L=2\pi$, $T=v\,\partial _x$ and $\mathcal L=\partial^2_v-v\,\partial_v$, in Fourier variables, we obtain $\lambda_m=\lambda_M=1$ and $C_M=(1+\sqrt3)/2$ according to~\cite[Section~II.1.3.2]{arnold2021}, so that $\lambda=1/(12+6\,\sqrt3)\approx0.0446582$. Using Fourier modes, a slightly better estimate is obtained from~\cite[Section~II.1.2]{arnold2021} with $\lambda\approx0.176048$.
\subsection{Direct spectral methods} 

In a series of papers, F.~Achleitner, A.~Arnold, E.~Carlen and several other collaborators use direct spectral methods. We refer in particular to~\cite{arnold2014sharp,achleitner2016linear,achleitner2017multi,achleitner2017optimal} and also~\cite{arnold2021} for an introduction to the method, which can be summarized as follows.

Let us consider~\eqref{Eqn:abstract} written after a Fourier transform in $x$, so that $T=i\,\xi\cdot v$, and acting on $\mathrm L^2(dx;\mathrm L^2_2)$ now considered as a space of complex valued functions. Assume that for some positive definite bounded Hermitian operator $P$ and some constant $\lambda\in(0,+\infty)$, we have
$$(L-T)^*P+P\,(L-T)\ge2\,\lambda\,P.$$
Let us consider the \emph{twisted} norm $\|\hat f\|_P^2 := \int_Q(\hat f,P\,\hat f)\,dx$ where $(\cdot,\cdot)$ is the natural extension of the scalar product as defined in~\eqref{scal}. From
$$\frac{d}{dt}\|\hat f\|_P^2 = -\,\langle \hat f, ((L-T)^*P + P(L-T))\hat f\rangle \leq -\,2\,\lambda\,\|\hat f\|_P^2,$$
for some $C>1$, we deduce that
$$C^{-1}\|f(t,\cdot,\cdot)\|_{\mathrm L^2(dx\,d\gamma_\alpha)}^2\le\|\hat f(t,\cdot,\cdot)\|_P^2 \le e^{-2\,\lambda\,t}\,\|\hat f_0\|_P^2\quad\forall\,t\ge0.$$
To our knowledge, $\mu$ has not yet been computed in the case of~\eqref{eq1}. The spectral decomposition 
$$h(t,x,v) = \sum_{\xi \in \mathbb Z^d} \sum_{k \in \mathbb{N}^d} a_{\xi,k}(t)H_k(v)\mathrm{e}^{-i\frac{2\pi}{L}\xi \cdot x}$$
provides an easy framework for finite dimensional approximations using the basis of Hermite functions $(H_k)_{k\in\N}$ and the numerical value $\mu\approx0.4$ has been obtained according to~\cite{Achleitner2020}.
\begin{table}[h]
    \centering
    \begin{tabular}{c|c}
    Proposition \ref{thm1vero}     & $\lambda \approx 0.07$  \\
    DMS \cite{dolbeault2015hypocoercivity}     & $\lambda \approx 0.04$ \\
    ADSW \cite{arnold2021} & $\lambda \approx 0.17$ \\
    Achleitner (numerics) \cite{Achleitner2020} & $\lambda \approx 0.4$ \\
    \end{tabular}
    \centering
    \caption{Comparison among different $\mathrm{L}^2$ hypocoercivity methods}
    \label{tab1}
\end{table}
\subsection{Comparison for decay rates in limit regimes}
Let $\alpha \geq 1.$
Corollary \ref{cor1} provides us with a decay estimate depending on the parameter $L,$ which represents the length of the spatial domain $\mathbb T.$ Note that \eqref{bestlionseq} is meaningful if $0<\tau<L.$ 
We shall now consider two situations, corresponding to $L \to \infty,$ where spatial diffusion dominates, and to $L \to 0,$ where the dominant term is the collision operator $\Delta_\alpha.$
In the first case, we have that the decay exponent 
$$\lambda \approx \frac{\tau}{L^3} \to 0, \qquad \text{ as } L \to \infty.$$
The hypocoercivity constant $C \approx 1.$ Hence, exponential decay is lost in the limit. 
On the other hand, for $L \to 0,$ we have
$$\lambda \approx \frac{1}{4 (P_\alpha +1) \, \|v_1 \, |v|^2\|^2_{\mathrm{L}^2_\alpha}} \approx 0.04,$$
if $\alpha =1.$
Moreover, $$C \approx 1.$$
This rate has the wrong order once compared to \cite{arnold2021}, where the authors recover the value $$\lambda \approx 1 - \sqrt{3/7}.$$
Our inaccuracy is mainly due to the incompatibility between \eqref{bestlionseq} and Lemma \ref{lemmaavg}. Moreover, the value of the Lions constant in \eqref{bestlionseq} is just an estimate and it is not expected to be as accurate as something achieved by a spectral method (even if its scaling is correct). 
 
\section*{Acknowledgements}
The author thanks F.~Achleitner and C.~Mouhot for stimulating discussions. This research project is funded by the European Union’s Horizon 2020 research and innovation program under the Marie Skłodowska-Curie grant agreement No 754362. Partial support has been obtained from the EFI ANR-17-CE40-0030 Project of the French National Research Agency. \\
{\small $\copyright$ 2022 by the author. Any reproduction for non-commercial purpose is authorized.}

\bibliography{biblio}{}
\end{document}